\documentclass[11pt]{article}
\usepackage{amsmath,amssymb,amsthm,amsfonts,graphicx,color,enumerate,float,epsfig}
\begin{document}

\newcommand{\kf}[1]{\marginpar{\tiny #1 --kf}}
\newcommand{\mb}[1]{\marginpar{\tiny #1 --mb}}
\newcommand{\kb}[1]{\marginpar{\tiny #1 --kb}}
\newtheorem{thm}{Theorem}[section]
\newtheorem{lemma}[thm]{Lemma}
\newtheorem{cor}[thm]{Corollary}
\newtheorem{prop}[thm]{Proposition}
\newtheorem{question}[thm]{Question}
\newtheorem*{mainthm0}{Theorem \ref{log}}

\newtheorem{remark}[thm]{Remark}
\newtheorem{important remark}[thm]{Important Remark}
\newtheorem{definition}[thm]{Definition}
\newtheorem{example}[thm]{Example}
\newtheorem{fact}[thm]{Fact}
\newtheorem{convention}[thm]{Convention}

\def\diam{\operatorname{diam}}
\def\vl{\operatorname{vl}}
\def\svl{\operatorname{svl}}
\def\cal{\mathcal}
\def\R{{\mathbb R}}
\def\Z{{\mathbb Z}}
\def\<{{\langle}}
\def\>{{\rangle}}

\def\Mod{{\rm Mod}}
\def\ep{{\varepsilon}}
\def\Isom{{\rm Isom}}
\def\C{{\mathcal C}}
\def\T{{\mathcal T}}
\def\TT{\overline{\mathcal T}}
\def\dTT{\partial \overline{\mathcal T}}
\def\S{{\mathcal S}}
\def\SS{{\overline{\mathcal S}}}
\def\M{{\mathcal M}}
\def\L{{\mathcal L}}
\def\D{\Delta}
\def\s{\sigma}
\def\t{\tau}
\newcommand{\cN}{{\cal N}}
\newcommand{\F}{{\mathbb F}}

\def\red{\textcolor{red}}

\title{The verbal width of acylindrically hyperbolic groups
is infinite}

\author{Mladen Bestvina, Ken Bromberg and Koji Fujiwara\thanks{This
    material is based upon work supported by the National Science
    Foundation under Grant No. DMS-1440140 while the authors were in
    residence at the Mathematical Sciences Research Institute in
    Berkeley, California, during the Fall 2016 semester.
The
    first two authors gratefully acknowledge the support by the National
    Science Foundation under grant numbers DMS-1607236 and DMS-1509171
    respectively. The third author is supported in part by
Grant-in-Aid for Scientific Research (No. 15H05739)}}

\maketitle

\begin{abstract}
We show that the verbal width is infinite for acylindrically hyperbolic
groups, which include hyperbolic groups, mapping class
groups and $Out(F_n)$.

\end{abstract}

\section{Verbal subgroups}
The {\em Brooks construction} is a method for constructing  essential quasi-morphisms on free groups. These Brooks quasi-morphisms can be extended to general acylindrically hyperbolic groups and in this note we use these quasi-morphisms to study {\em verbal subgroups}. A verbal subgroup generalizes the notion of a commutator subgroup. We begin with a precise definition. Suppose that $G$ is a group and $\mathbb F_k$ is the free group
generated by $x_1,\cdots,x_k$. For any choice of $w\in \mathbb F_k$
substitution defines a map, also denoted $w$:
$$w:G^k\to G.$$
In what follows we always assume that $w$ is a non-trivial.

The image of this function is the {\it verbal subset}, and is denoted
by $w[G]$.  The {\it verbal subgroup} $w(G)<G$ is the subgroup
generated by $w[G]$. For example, if $k=2$ and
$w=x_1x_2x_1^{-1}x_2^{-1}$ then $w[G]$ is the set of commutators in
$G$ and $w(G)=[G,G]$ is the commutator subgroup of $G$.

Let $e_i$ be the sum of the exponents of $x_i$ in $w$.
For example, if $w=x_2 x_1 x_2^{-2}$, then $e_1=1, e_2=1-2=-1,
e_3= \cdots =e_k=0.$
Let $d(w) \ge 1$ be the g.c.d. of $e_i$'s.
If they are all $0$, define $d(w)=0$. Note that the condition that $d(w) = 0$ implies that $w(G) \subset [G,G]$ for when $d(w)=0$, $w(G)$ will be in the kernel of any homorphism from $G$ to an abelian group.

For each $g \in w(G)$, 
define its {\it verbal length} by
$$vl_w(g)= \min \{n | g=g_1 \cdots g_n, g_i \in w[G]^{\pm1} \}.$$ The
{\it width} of $w(G)$ is the supremum of $vl_w(g)$ over all $g \in w(G)$.
Note that if $d(w)>1$, then $g^d \in w[G]$. (Replace each $x_i$ in $w$ with $g^{a_i e_i}$ where $\sum a_i e_i = d$.) In particular if $d =1$ then $w[G] = w(G) = G$ and the width is $1$.
The reader can consult the book \cite{segal} for more information on
the subject.  If $w$ is the commutator $[x_1,x_2]$, the verbal length
is called the {\it commutator length}.

A group $G$ is {\em acylindrically hyperbolic} if it has a
non-elementary acylindrical action on a $\delta$-hyperbolic
space \cite{osin-acylindrical}. Recall that an isometric action of $G$
on a metric space $X$ is {\em acylindrical} if for all $D>0$ there exist
$L,N>0$ such that if $d(x,y)>L$ then the set
$$\{g \in G | d(x, gx) < D \mbox{ and } d(y, gy)<D \}$$ has $<N$
elements. The first nontrivial example is due to Bowditch
\cite{bhb:tight} who showed that the action of the mapping class
group on the curve complex is acylindrical. There are now many
examples with the key point being that many seemingly weaker geometric
criteria imply that the group is acylindrically hyperbolic. See
\cite{osin-acylindrical} and \cite{bbf-bc}.

Here is our main  result. 
\begin{thm}\label{width}
Suppose that $G$ is acylindrically hyperbolic and that $d(w) \neq 1$. 
Then the width of $w(G)$ is infinite.
\end{thm}
This result generalizes the work of Rhemtulla \cite{rhem} and
Myasnikov-Nikolaev \cite{nm} who proved the theorem for free groups
and hyperbolic groups, respectively.

Similarly to the stable commutator length, 
one can define the {\it stable verbal length} of $g \in w(G)$, $svl_w(g)$,
as follows:
$$svl_w(g) = \liminf _{n \to \infty} \frac{vl_w(g^n)}{n}.$$
If $svl_w(g)>0$ for some $g$ then $w(G)$ have infinite verbal width. However, if $d(w) \ge 1$ this method cannot be used due to the following lemma:
\begin{lemma}[Calegari-Zhuang \cite{cz}]
If $d=d(w) \ge 1$ then $vl_w(g^n)$ is bounded and $svl_w(g) = 0$ for all $g \in w(G)$.
\end{lemma}

We may supress $w$ and write $vl, svl$ instead of $vl_w, svl_w$.
\proof
As observed above, $g^d$ is in $w[G]$ for all $g \in G$ so $vl(g^{nd}) = 1$. Since $vl(gh) \le vl(g) + vl(h)$ this implies that $vl(g^n)$ is bounded and $svl(g) = 0$.
\qed

On the other hand, when $d(w) =0$ we have the following which implies Theorem \ref{width} in this case.
\begin{thm}\label{svl}
If $G$ is acylindrically hyperbolic and $d(w)=0$, then $svl_w(g) >0$ for some element $g \in w(G)$.
\end{thm}
If $G$ is a free group then this is Corollary 2.16 of \cite{cz}.

\subsection{An outline of proof of Theorem 1.1 for free groups}

To illustrate the main idea, we sketch the proof of Theorem 1.1 in the
case that $G=F$ is the free group with basis $\{a,b\}$ and
$w=x_1x_2x_1x_2^{-1}$, so that $d=d(w)=2$. Consider
$g_i=ab^{2i}ab^{2i+1}\in G$ for $i=1,2,\cdots$. This sequence has the
property that distinct occurrences of any $g_i$ in any reduced word
have trivial overlap. Denote by $H_i:G\to\Z$ the Brooks counting
quasi-morphism with respect to $g_i$. For any $y,z\in G$ we have
$|H_i(yz)-H_i(y)-H_i(z)|\leq 3$ by the usual tripod argument, since at
most 3 copies of $g_i$ along a tripod can have the tripod point in the
interior. The key observation now is that in fact
$H_i(yz)-H_i(y)-H_i(z)=0$ for all but at most 3 values of $i$ by the
non-overlapping property of the $g_i$'s (the exceptional values of
$i$ depend on $y$ and $z$).

Now suppose that $g\in G$ has $vl_w(g)=1$, so
$g=x_1x_2x_1x_2^{-1}$ for some $x_1,x_2\in G$. Then
$$H_i(g)=H_i(x_1)+H_i(x_2)+H_i(x_1)+H_1(x_2^{-1})=2H_i(x_1)$$
is even for all but $3\times 3=9$ values of $i$. Thus to detect $g$
with $vl_w(g)>1$ it suffices to ensure that $H_i(g)$ is odd for 10
values of $i$. Similarly, to detect that $vl_w(g)$ is large it
suffices to ensure that $H_i(g)$ is odd for sufficiently many $i$. An
element such as
$$(ab)^2(ab^2)^2(ab^3)^2(ab^4)^2\cdots (ab^{N-1})^2(ab^N)^2$$
will do.

For a general acylindrically hyperbolic group $G$ we perform the above
construction on a suitable Schottky subgroup $F\subset G$ and then use
the method of Hull and Osin \cite{ho} to extend the quasi-morphism
from $F$ to $G$. We will review their construction in Section 4 and show
that the key observation above continues to hold for the extended
quasi-morphisms. 

\section{Extending Brooks quasi-morphisms to acylindrically hyperbolic groups}
\label{section.review}
We first recall the definition of a quasi-morphism. Let $G$ be a group. Then
$$H:G \to \R$$
is a {\em quasi-morphism} if 
$$\underset{\alpha, \beta \in G}{\sup} |H(\alpha\beta) - H(\alpha) - H(\beta)| = \Delta(H) < \infty.$$
The constant $\Delta = \Delta(H)$ is the {\em defect} of $H$.
Note that if $\Delta =0$ then $H$ is a homomorphism. A quasi-morphism is {\em anti-symmetric} if $H(-\alpha) = -H(\alpha)$.

One way to construct a quasi-morphism that is not a homomorphism is to start with a homomorphism and then add on a bounded function. Of course, this is not an interesting example. The {\em Brooks construction} is a way of building an anti-symmetric quasi-morphisms that are not a bounded distance from a homomorphism.

Let $F = \langle a,b \rangle$ be the free group on two generators and
let $w$ be a reduced word in $F$. 
For $x \in F$ let $N_w(x)$ be the
number of copies of $w$ in $x$ when $x$ is written as a cylically reduced word
and let $H_w(x) = N_w(x) - N_{w^{-1}}(x)$. Note that $H_w(-x) = -H_w(x)$ and $H_w(w^n) = n$. Brooks proved the
following:
\begin{thm}\label{brooks}
The function $H_w$ is an anti-symmetric quasi-morphism and if $w$ is
not a power of $a$ or $b$ then it is not a bounded distance from a
homomorphism. If $w$ is cyclically reduced, then $H_w(w^n)\geq n$.
\end{thm} 
Note that if $N$ is a finite group then the Brooks quasi-morphisms can be extended to $F \times N$ by choosing them to be constant on the second factor.

In \cite{dgo}, Dahmani-Guirardel-Osin show that an acylindrically
hyperbolic group $G$ contains a copy of a {\em hyperbolically
  embedded} $F \times N$ where $N$ is the maximal finite normal
subgroup of $G$. In \cite{ho}, Hull-Osin show that any anti-symmetric
quasi-morphism on a hyperbolically embedded subgroup extends to a
quasi-morphism of the entire group. Combing these two results we have
the following theorem.

\begin{thm}\label{quasi.tree}
Let $G$ be acylindrically hyperbolic. Then there exists a free group $F = \langle a,b\rangle < G$ such that for every Brooks quasi-morphism $H_w$ there is a quasi-morphism $H: G \to \R$ such that $H|_F = H_w$.
\end{thm}

\begin{remark}
There is a weaker version of this theorem (that would be good enough for our applications here) that follows from \cite{bbf-bc}. The approach in \cite{bbf-bc} is more direct as it does not go through the theory of hyperbolically embedded subgroups. We also note that both approaches use the projection complex from \cite{bbf} in an essential way.
\end{remark}

As a demonstration of our methods we first give a proof of Theorem \ref{svl}.
\proof
Given a quasi-morphism $H$ with the defect $\Delta$ and an element $g = g_1\dots g_n$ by repeatedly applying the quasi-morphism bound we have
$$\left|H(g) - \sum_{i=1}^n H(g_i)\right| \le (n-1)\Delta.$$
If $g= w(g_1, \dots, g_k)$ and $H$ is anti-symmetric this becomes
$$\left|H(g) - \sum_{i=1}^n e_i H(g_i)\right| \le (|w|-1)\Delta$$ and
when $d(w) =0$ (so all the $e_i=0$) this becomes $|H(g)| \le (|w|-1)
\Delta$ for $g \in w[G]$. More generally for $g \in w(G)$ we have
$|H(g)| \le (vl(g)|w|-1) \Delta$ and therefore if $|H(g)|>0$ we have
$svl(g)>0$ since $H(g^n) \ge n H(g)$ for all $n>0$.

We will use the Brooks construction (and the Hull-Osin extension) to
find a $g \in w(G)$ with $H_g(g^n) \geq n$. To do this we need to find
a cyclically reduced word in $w(F) \subset w(G) \cap F$.
Pick a non-trivial element $h \in w(F)$. If it is cyclically reduced
let $g=h$ and we are done. If not, then $h=a\cdots a^{-1}$ for a basis
element $a$ (or its inverse). Let $h'$ be obtained from $h$ by
swapping $a$'s and $b$'s (with $b$ another basis element). Then $hh'$ is cyclically reduced and still in $w(F)$ so $g =hh'$ is the desired element. \qed

If $w \in [G,G]$ then for any $g \in w(G)$, $cl(g) \le cl(w) vl(g)$ so $scl(g) \le cl(w) svl(g)$. In particular if $scl(g)>0$ then $svl(g)>0$ and Theorem \ref{svl} would follow if we knew that every verbal subgroup of an acylindrically hyperbolic group had an element $g$ with $scl(g)>0$. However, proving this does not seem any easier that the more general proof above.

One can also ask if $scl(g)=0$ implies that $svl_w(g) = 0$ for all
$w$. Here the answer is negative. For example, take $w=[[x,y],[z,u]]$
and
$$G=\langle a,b,c,d,t\mid t[[a,b],[c,d]]t^{-1}=[[a,b],[c,d]]^{-1}\rangle$$
Then for $g=[[a,b],[c,d]]$ we have that $g$ is conjugate to $g^{-1}$,
which forces $scl(g)=0$. On the other hand, we claim that
$svl_w(g)>0$. Indeed, if $svl_w(g)=0$, then $g^n$ can be written as a
product of a sublinear number of double commutators, which would imply
that $scl_H(g)=0$ where $H=[G,G]$. We now argue that $scl_H(g)>0$. In
fact we will show that there is an index 2 subgroup $N<G$ with $H<N$
and so that $N$ surjects to the free group $F_4=\<a,b,c,d\>$ with $g$
mapping to $[[a,b],[c,d]]$. Since nontrivial elements in free groups
have positive $scl$ the claim follows. In fact, it shows
$scl_N(g)>0$. It immediately implies $scl_H(g)>0$.

We let $N$ be the kernel of $G\to\Z/2$ that sends $t$ to 1 and
$a,b,c,d$ to 0. The corresponding double cover $Y$ of the presentation
2-complex $X$ of $G$ consists of the disjoint union of two roses $R_i$
with petals labeled $a_i,b_i,c_i,d_i$, $i=1,2$, with edges $t_i$
connecting the vertex of $R_i$ to the vertex of $R_{3-i}$. The map to
$X$ is the obvious one, sending $a_i$ to $a$ etc. The relation 2-cell
in $X$ 
lifts to two 2-cells in $Y$, with attaching maps
$t_i[[a_{i+1},b_{i+1}],[c_{i+1},d_{i+1}]]t_i^{-1}[[a_i,b_i][c_i,d_i]]$
with indices taken mod 2. Now map the $1$-skeleton of $Y$ to the rose corresponding to
$\<a,b,c,d\>$ via $t_i\mapsto *$, $a_1,c_2 \mapsto a$,
$b_1,d_2 \mapsto b$, $c_1,a_2 \mapsto c$, $d_1,b_2 \mapsto
d$. We then extend this to the two $2$-cells.
This is possible since via the attaching maps the boundary of 
the $2$-cells are mapped
to $[[c,d],[a,b]][[a,b],[c,d]]$ and $[[a,b],[c,d]][[c,d],[a,b]]$,
which are trivial in $\<a,b,c,d\>$ since $[x,y]^{-1}=[y,x]$.

\section{Some facts about the Hull-Osin extension}
Unfortunately, rather than just the statement of Theorem \ref{quasi.tree} we need some elements of the proof in \cite{ho}. We review them now. In this section $F$ can be any hyperbolically embedded subgroup in $G$.

It is convenient to replace the quasi-morphism with a function 
on $G \times G$, called a bicombing in \cite{ho}.
If $H$ is a quasi-morphism we define $r(x,y) = H(x^{-1}y)$. Note that $r(zx,zy) = r(x,y)$ and $|r(x,y) + r(y,z) + r(z,x)|$ is bounded by the defect of $H$. On the other hand if we are given a map $r(x,y)$ (satisfying the properties from the previous sentence) then the map $x \mapsto r(1,x)$ is a quasi-morphism so $r$ determines $H$ just as $H$ determines $r$. In particular, to construct $\tilde H$, in \cite{ho} they first construct $\tilde r\colon G \times G \to \R$. To construct $\tilde r$ for each $x,y \in G$  and each coset $aF$ is associated a finite collection of pairs $E(x,y; aF) = \{(x_i, y_i)\}$ where $x_i,y_i \in F$. For the convenience of the reviewer we briefly review the construction of the sets and then  state the key properties that we will need.

Let $\Gamma$ be a Cayley graph for $G$ formed from a generating set that contains every element of $F$. Given $x,y \in G$ let $\gamma$ be a geodesic in $\Gamma$ from $x$ to $y$. Each $F$-coset has diameter one in $\Gamma$ so $\gamma$ will intersect a given coset $aF$ in at most two points $x'$ and $y'$. We say that $\gamma$ {\em essentially penetrates} $aF$ if any path in $\Gamma$ from $x'$ to $y'$ that doesn't contain any $F$-edges has length $\ge C$ where $C$ is a constant that only depends on $G$ and $F$. We let $S(x,y)$ be the set of cosets $aF$ where there is some geodesic from $x$ to $y$ that essentially penetrates $aF$. A central fact from \cite{ho} is that if there is one geodesic that essentially penetrates then every geodesic from $x$ to $y$ must intersect $aF$. For each coset in $aF\in S(x,y)$ we let $E(x,y; aF)$ be the set of pairs $(x',y') \in F$ such that $ax'$ and $ay'$ are the entry and exit points for a geodesic from $x$ to $y$ in $\Gamma$. For each coset the particular representative $a$ is not important except that the choice needs to be fixed for once and all. If $aF\not\in S(x,y)$ then $E(x,y;aF)$ is empty.

We now define
$$\tilde r(x,y) = \sum_{aF\in S(x,y)}\left(\frac1{|E(x,y; aF)|} \sum_{(x',y') \in E(x,y; aF)} r(x',y')\right).$$
For this to be well defined we need the sum to be finite. The inside sum is finite by Lemma 3.8 of \cite{ho} and the outside sum is finite since $S(x,y)$ is finite by Corollary 3.4. Note that while in \cite{ho} it is only stated that the size of $E(x,y;aF)$ is finite it is in fact uniformly bounded which will be important in the proof of Lemma  
\ref{bad cosets} later.

The following lemma is a combination of Lemma 3.9 and (the proof of) Lemma 4.7 in \cite{ho}.
We fix the word metric with respect to a finite 
generating set  on $F$ and denote the distance
between $x,y$ by $|x-y|$.
\begin{lemma}\label{technical fact}
Given $x,y,z \in G$ for all but at most two cosets $aF$ exactly one of the following three possibilities holds:
\begin{enumerate}
\item $E(x,y; aF) = \emptyset$;

\item $E(x,y; aF) = E(x,z; aF)\not= \emptyset $ and $E(y,z; aF) =\emptyset$;

\item $E(x,y; aF) = E(y,z; aF)\not= \emptyset$ and $E(x,z; aF) = \emptyset$.
\end{enumerate}
If $aF$ doesn't satisfy the above then either
\begin{enumerate}[(A)]
\item All of $E(x,y; aF)$, $E(x,z; aF)$ and $E(y,z; aF)$ are non-empty and for any pairs $(x', y') \in E(x,y;aF)$, $(x'', z') \in E(x,z; aF)$ and $(y'', z'') \in E(y,z; aF)$, $|x'- x''|$, $|y'- y''|$ and $|z'-z''|$ are uniformly bounded.

\item Only $E(x,z; aF)$ is empty and for any pairs $(x', y') \in E(x,y;aF)$ and $(z',y'') \in E(z,y;aF)$, $|x'-z'|$ and $|y'-y''|$ are uniformly bounded or the same statement holds with $y$ and $z$ swapped.

\item Only $E(x,y; aF)$ is non-empty and for all pairs $(x',y') \in E(x,y;aF)$, $|x'-y'|$ is uniformly bounded.
\end{enumerate}
\end{lemma}

Given a pair $(x,y) \in G$ let $B(x,y)$ be the collection of cosets that don't satisfy (1)-(3) and let $B(x,y,z)$ be the union of $B(x,y)$, $B(y,z)$ and $B(z,x)$. By Lemma \ref{technical fact}, $B(u,v)$ contains at most 2 cosets so $B(x,y,z)$ contains at most six. Cosets in $B(x,y,z)$ are of type (A), (B) or (C) depending on which of the conditions in Lemma \ref{technical fact} they satisfy.

We are interested in the sum $\tilde r(x,y) + \tilde r(y,z) + \tilde r(z,x)$. It will be convenient to define new  sets $E(x,y,z;aF)$ to be the product of the sets of pairs $E(x,y;aF)$, $E(y,z;aF)$ and $E(z,x;aF)$. Note that one or more of the sets may be empty in which case the product would be empty. (In fact for at most one coset at least one of the sets will be empty.) To get around this if $E(u,v;aF)$ is empty we make it non-empty by adding the ``empty pair'' $(\emptyset, \emptyset)$ and we define $r(\emptyset, \emptyset) =0$. With this modification $E(x,y,z;aF)$ will always be a triple of pairs in $F \cup \{\emptyset\}$. Next we define
$$\rho(x,y,z;aF) = \frac{1}{|E(x,y,z;aF)|}\sum_{E(x,y,z;aF)} r(x_-, y_+) + r(y_-,z_+) + r(z_-,x_+)$$
and observe that
$$\tilde r(x,y) + \tilde r(y,z) + \tilde r(z,x) = \sum_{aF} \rho(x,y,z;aF).$$
To show that $\tilde r$ determines a quasi-morphism Hull-Osin show that for nearly all cosets the expression $\rho(x,y,z;aF)$ is zero and for the finitely many when it is not it is uniformly bounded.

\begin{cor}\label{zero}
If $aF \not\in B(x,y,z)$ then $\rho(x,y,z;aF) = 0$. 
\end{cor}

\proof
If $aF \not\in B(x,y,z)$ then either $E(x,y,z;aF)$ is the triple of empty pairs and $\rho(x,y,z;aF) = 0$ or all the terms in the sum cancel and again $\rho(x,y,z;aF) = 0$. 
\qed

\section{Many independent quasi-morphisms}
For the remainder of the paper we can assume that $d(w)>1$. If $H: G \to \Z$ is a homomorphism, an easy calculation gives that for any $g \in w(G)$ we have that $H(g)$ is divisible by $d(w)$.  We will construct a family of quasi-morphisms where this is true for nearly all the quasi-morphisms in the family where the number of exceptions is bounded above by the $vl(g)$.

Let $F\times N$ be hyperbolically embedded
in $G$ where $F$ is the free group of rank at least two
and $N$ is a finite group. 
For simplicity we suppose the rank of $F$ is two
in the following.
We now fix a sequence of words that we will use to build  Brooks' quasi-morphisms on $F$, then 
extend it to $F \times N$, trivially on $N$.
 Let $g'_i = ab^{2i}$, $g''_i = ab^{2i+1}$ and $g_i = g'_i g''_i$ and let $H_i = H_{g_i}$ be the Brooks quasi-morphism and $r_i$ the corresponding bicombings.
 
 We fix the word metric with respect a finite
 generating set  on $F \times N$ and denote
 the distance between $x,y$ by $|x-y|$.
 
\begin{lemma}\label{bad cosets}
Given a triple of pairs $(x_-, x_+), (y_-,y_+), (z_-,z_+)$ in $F\times N$ with $|x_- - x_+|$, $|y_- - y_+|$ and $|z_- - z_+|$ bounded by $L$ there are at most
\begin{itemize}
\item $L$ of the $r_i$ such that $r_i(x_-,x_+) \neq 0$;

\item $2L$ of the $r_i$ such that $r_i(x_-,y_+) + r_i(y_-,x_+) \neq 0$;

\item  $3L+3$ of the $r_i$ such $r_i(x_-,y_+) + r_i(y_-, z_+) + r_i(z_-,x_+) \neq 0$.
\end{itemize}
Therefore there is a uniform bound on the number of $\rho_i$ where $\rho_i(x,y,z;aF) \neq 0$.
\end{lemma}

\proof
We only discuss the case that $N$ is trivial.
The general case is similar.

If $r_i(x_-,x_+) \neq 0$ there is a translate of the word $g_i$ in the segment between $x_-$ and $x_+$ in the Cayley graph (with the standard generators). Since two $g_i$ can't intersect in a segment (a very bad) upper bound for the number of $r_i$ with $r_i(x_-,x_+)\neq 0$ is $|x_- - x_+|\le L$.

The triple $x_-$, $x_+$ and $y_+$ form a tripod in the Cayley graph and let $m$ be the central vertex. Then $r_i(x_-, y_+) + r_i(y_+,x_+) = 0$ unless there is a a translate of $g_i$ in the segment from $y_+$ to $x_-$ that intersects the segment from $m$ to $x_-$ or a translate in the segment from $y_+$ to $x_+$ that intersects the segment from $m$ to $x_+$. Again using the fact that two $g_i$'s can't overlap in a segment an upper bound for the number of $r_i$'s with $r_i(x_-, y_+) + r_i(y_+,x_+) \neq 0$ is $|x_- - m|+ |x_+ - m| = |x_- - x_+| \le L$. Similarly there are at most $L$ of the $r_i$ such that $r_i(x_+, y_+) + r_i(y_+,x_+) \neq 0$ or equivalently $r_i(y_+, x_+)= r_i(y_-,x_+)$ for all but $L$ of the $r_i$ and therefore $r_i(x_-,y_+) + r_i(y_-,x_+) = 0$ for all but $2L$ of the $r_i$.

Now we examine the tripod formed by $x_-, y_-$ and $z_-$. As with the original Brooks' argument the sum
$$r_i(x_-,y_-) + r_i(y_-, z_-) +r_i(z_-, x_-)$$
is zero unless the a translate of the word $g_i$ intersects the central vertex of the tripod. At most three such words can intersect the central vertex so the sum is non-zero for at most 3 of the $r_i$. As above for at most $L$ of the $r_i$ we have $r_i(x_-, y_-) = r_i(x_i,y_+)$, etc. Therefore
$$r_i(x_-,y_+) + r_i(y_-, z_+) +r_i(z_-, x_+)=0$$
for all but at most $3L+3$ of the $r_i$.

Since $F$ is hyperbolically embedded in $G$,
if $aF \not\in B(x,y,z)$ then $\rho_i(x,y,z;aF)=0$ for all $i$ by Corollary \ref{zero}. For cosets $aF \in B(x,y,z)$ of type (A) there will be at most $(3L+3)|E(x,y,z;aF)|$ of the $\rho_i$ with $\rho_i(x,y,z;aF) \neq 0$, for cosets of type (B) at most $2L|E(x,y,z;aF)|$ and for cosets of type (C) at most $L|E(x,y,z;aF)|$. By Lemma 3.8 of \cite{ho} $|E(u,v;aF)|$ is uniformly bounded\footnote{In \cite[Lemma 3.8]{ho} it is only claimed that $|E(x,y;aF)|$ is finite however it is easy to see that their proof shows that the bound is uniform since the constant $C$ in Lemma 3.3/2.4 is uniform.} and therefore so is $|E(x,y,z;aF)|$. It follows that there is uniform bound on the number of $\rho_i$ with $\rho_i(x,y,z;aF)\neq 0$.\qed

Since $F \times N$ is hyperbolically embedded in $G$,
let $\tilde H_i:G \to \R$ be the Hull-Osin extension of the $H_i$ and $\tilde r_i$ the corresponding bicombings.

\begin{prop}\label{almost.homomorphism}

There exists an $M>0$ such that for
any $x,y \in G$, 
$$\tilde H_i (xy)-\tilde H_i(x)-\tilde H_i(y) =0$$
holds except for at most $M$ of the  $\tilde H_i$.
It follows that for any $a_1, \cdots a_k \in G$, 
$$\tilde H_i(a_1 \cdots a_k) -\tilde H_i(a_1) - \cdots -  \tilde H_i(a_k) =0$$
holds except for at most $M(k-1)$ of the $\tilde H_i$.

\end{prop}

\proof
First observe
$$\tilde H_i (xy)-\tilde H_i(x)-\tilde H_i(y) = \tilde r_i(id, xy)  + \tilde r_i(x, id)
+ \tilde r_i(xy, x)$$
so we can instead show that
$$\tilde r_i(x,y) + \tilde r_i(y,z) + \tilde r_i(z,x) = 0$$
for all but $M$ of the $\tilde r_i$. But this follows from Lemma \ref{bad cosets}, as for all but at most 6 cosets $\rho_i(x,y,z;aF) = 0$ for all $i$ and for each of these bad cosets there is a uniform bound on the number of $\rho_i$ with $\rho_i(x,y,z;aF) \neq 0$.
\qed
\\

\begin{lemma}\label{special.word}
For each $K$ there exists $g \in w(F\times N)$
such that 
$H_i(g)=1$ for all $1 \le i \le K$.
\end{lemma}

\proof 
We will find an element $g \in w(F)$.
Recall that we are assuming that $d=d(w)>1$ and therefore for any $f \in F$, we have $f^d \in w[F]$. Let $h'_i = (g'_i)^d$ and $h''_i = (g''_i)^d$. Then the product $h_i = h'_i h''_i$ contains a single copy of $g_i$. Let $g = h_1 h_2 \cdots h_K$. Note that $g$ is already reduced since there are only positive powers of $a$ and $b$ in the $h'_i$ and $h''_i$. Furthermore by our construction of the $g_i$ there will be exactly one copy of $g_i$ in $g$ and no copies of $g^{-1}_i$. Therefore $H_i(g) = 1$ for $1\le i \le K$. \qed
\\

We now give a proof of Theorem \ref{width} when $d(w)>1$.

\proof
We first show that for any $g \in w(G)$ we have that $\tilde H_i(g)$ is divisible by $d(w)>1$ for all but a bounded number of the $\tilde H_i$ where the bound only depends on $vl(g)$. To see this we first observe that if $g = w(g_1,\dots, g_k) \in w[G]$ then by Proposition \ref{almost.homomorphism}
$$\tilde H_i(w(g_1, \cdots, g_k)) = \sum _{j=1}^k e_j \tilde H_i(g_j),$$
for all but $M(|w|-1)$ of the $\tilde H_i$. In particular, for any $g \in w[G]$ there are at most $M(|w|-1)$ of the $i$ such that $\tilde H_i(g)$ isn't divisible by $d(w)$. Similarly if $g \in w(G)$ is product of $vl(g)$ elements $g_j \in w[G]$ then
$$\tilde H_i(g) = \sum_{j=1}^{vl(g)} \tilde H_i(g_j)$$
for all but $M(vl(g) - 1)$ of the $\tilde H_i$. If all of the $\tilde H_i(g_j)$ are divisible by $d(w)$ then so is $\tilde H_i(g)$ so we have that $\tilde H_i(g)$ is divisible by $d(w)$ for all but at most $M(vl(g) - 1) + vl(g)(M(|w|-1))$ of the $\tilde H_i$. In particular a bound on $vl(g)$ gives a bound on the number of $\tilde H_i$ where $\tilde H_i(g)$ is not divisible by $d(w)>1$.

On the other hand, by Lemma \ref{special.word}, for any $K>0$ we can find a word $h_K$ such that $\tilde H_i(h_K) = H_i(h_K)=1$ for $1\le i \le K$. Therefore $vl(h_K)\to \infty$ as $K\to \infty$.
\qed
\\
From the above proof we see that $vl(h_K) \ge K/(M|w|)$. We know of no examples where this bound is sharp.

\bibliography{./ref2.bib}
\end{document}